\documentclass[10pt]{article}               
\usepackage{amsmath,amsthm,amssymb}                 

\theoremstyle{plain}
\newtheorem{theorem}{Theorem}[section]

\newtheorem{lemma}[theorem]{Lemma}
\newtheorem{proposition}[theorem]{Proposition}
\newtheorem{fact}{Fact}[section]

\theoremstyle{definition}
\newtheorem{definition}[theorem]{Definition}
\newtheorem{example}[theorem]{Example}

\numberwithin{equation}{section}

\def\Bbb#1{{\fam\msbfam\relax#1}}
\font\fivemsb=msbm5
\font\sevenmsb=msbm7
\font\tenmsb=msbm10
\newfam\msbfam
\textfont\msbfam=\tenmsb
\scriptfont\msbfam\sevenmsb
\scriptscriptfont\msbfam\fivemsb

\def\Bbb{\rm \bold}
\def\spc{{\Bbb C}}
\def\spr{{\Bbb R}}

\def\vf{\varphi}

\def\a{\alpha}
\def\b{\beta}

\def\vF{\varPhi}

\def\k{\kappa}
\def\l{\lambda}
\def\m{\mu}
\def\n{\nu}
\def\om{\omega}
\def\r{\rho}

\def\t{\theta}

\def\x{\xi}
\def\vX{\varXi}
\def\y{\eta}

\def\cd{{\cal D}}
\def\ce{{\cal E}}

\def\cl{{\cal L}}

\def\cs{{\cal S}}
\def\cw{{\cal W}}
\def\la{\langle}
\def\ra{\rangle}
\def\lla{\langle\!\langle}
\def\rra{\rangle\!\rangle}

\begin{document}

\begin{center}
{\LARGE{\bf Characterization of 
Product Measures\\
\medskip by Integrability 
Condition}}
\end{center}

\medskip
\begin{center}
{\rm Nobuhiro ASAI\\
International Institute for Advanced Studies\\
Kizu, Kyoto
619-0225, Japan.}\\
{\tt asai@iias.or.jp}

\end{center}
\medskip



\section{Introduction}
Let $(\ce',\m)$ be the real Gaussian space, 
where $\ce'$ is the space of tempered distributions
and $\m$ be the standard Gaussian measure on $\ce'$.
In the recent papers \cite{akk3,akk6,akk7}, 
Asai, Kubo and Kuo (AKK for short) have shown that in order to construct 
the Gel'fand triple $[\ce]_u\subset L^2(\ce',\m)\subset [\ce]^{*}_u$
associated with a growth function $u\in C_{+,1/2}$, 
essential conditions on $u$ are (U0)(U2)(U3) 
stated in Section \ref{sec:2-2}.
Legendre transform and dual Legendre transform (Section \ref{sec:2-1})
play important roles to get this result.
We note that Gannoun et al. \cite{ghor}
have obtained similar results independently. 
Some relationships with \cite{ghor} are discussed in Section \ref{sec:2-2}.
In addition, the intrinsic topology for $[\ce]_u$
has been given and 
the characterization theorem for positive Radon measures
on $\ce'$
has also been proved by considering an integrability 
condition \cite{akk4,akk7}.

Now it is natural to ask whether ``positivity" of white
noise operators can be discussed in some sense and characterized.
To answer this question, we consider the Gel'fand triple 
over the Complex Gaussian space $(\ce'_c,\m_c)$,
i.e. $\ce'_c=\ce'+i\ce'$ equipped with the product measure
$\m_c=\m'\times\m'$ where $\m'$ is the Gaussian measure
on $\ce'$ with variance $1/2$ (Section \ref{sec:2-2}). 
Following AKK's Legendre transform technique, 
we have 
$\cw_{u_1,u_2}\subset L^2(\ce'_c,\m_c)\subset [\cw]^{*}_{u_1,u_2}$
for functions $u_1,u_2\in C_{+,1/2}$ 
satisfying (U0)(U2)(U3).
Several examples for $u_1, u_2$ are given in Section \ref{sec:2-3}.
We remark that Ouerdiane \cite{oue} studied a special case 
$u_1(r^2)=u_2(r^2)=\exp(k^{-1}r^k)$, where $1\leq k\leq 2$.
In Section \ref{sec:3}, the characterization theorem for measures 
can be extended to the case of positive 
product Radon measures on $\ce'\times \ce'$.  
In addition, the notion of 
pseudo-positive operators is naturally introduced 
via kernel theorem and characterized by an integrability
condition.  Lemma \ref{lem:3-2} plays crucial roles in 
Section \ref{sec:3}.
 
\section{White Noise Functions}
\subsection{Legendre Transform and Dual Legendre Transform}
\label{sec:2-1}
In this section we introduce the Legendre transform and 
dual Legendre transform which will be
used for the constructions of the Gel'fand triples 
over the real and complex Gaussian spaces. 

First, let us define two kinds of convex functions.
A positive function $f$ on $[0, \infty)$ 
is called
\begin{itemize}
\item[(1)] {\em (log, exp)-convex} if the function 
 $\log f(e^{x})$ is convex on $\spr$;
\item[(2)] {\em (log, $x^{k}$)-convex} if the function 
 $\log f(x^{k})$ is convex on $[0, \infty)$. Here $k>0$.
\end{itemize}   

\smallskip
Let $C_{+, \log}$ denote the collection of all positive 
continuous functions $u$ on $[0, \infty)$ satisfying the 
condition:
\begin{equation}
 \lim_{r\to\infty} {\log u(r) \over \log r}=\infty. \notag
\end{equation}
The {\em Legendre transform} $\,\ell_{u}$ of $u \in 
C_{+, \log}$ is defined to be the function
\begin{equation}
 \ell_{u}(t) = \inf_{r>0} {u(r) \over r^{t}}, \qquad
  t\in [0, \infty). \notag
\end{equation}

\smallskip
 Let $C_{+, \,1/2}$ denote the collection of all positive 
 continuous functions $u$ on $[0, \infty)$ satisfying the 
 condition:
 \begin{equation}
  \lim_{r\to\infty} {\log u(r) \over \sqrt{r}}=\infty. \notag
 \end{equation}
The {\em dual Legendre transform} $\,u^{*}$ of $u \in 
C_{+, 1/2}$ is defined to be the function
\begin{equation}
 u^{*}(r) = \sup_{s\geq 0} {e^{2\sqrt{rs}} \over u(s)}, 
   \qquad  r\in [0, \infty). \notag
\end{equation}
Note that $C_{+, 1/2} \subset C_{+, \log}$. 
Assume that $u\in C_{+, \log}$ and  
$\lim_{n\to\infty} \ell_{u}(n)^{1/n} =0$. We define the 
$L$-{\em function} $\cl_{u}$ of $u$ by
\begin{equation} \label{eq:2-1}
  \cl_{u}(r) = \sum_{n=0}^{\infty} \ell_{u}(n) r^{n}.
\end{equation}

For discussions in the rest of the paper, we will
need the following facts in \cite{akk6,akk7}.
See also \cite{asai1}.

\begin{fact} \label{fact:2-1} 
 (1) Let $u\in C_{+, \log}$ be (log, exp)-convex. Then
 its $L$-function $\cl_{u}$ is also (log, exp)-convex and 
 for any $a>1$,
 \begin{equation} 
  \cl_{u}(r) \leq {ea \over \log a} u(ar), \qquad
  \forall r\geq 0.  \notag
 \end{equation}
 \par\noindent
 (2) Let $u\in C_{+, \log}$ be increasing and 
 (log, $x^{k}$)-convex. Then there exists a constant $C$, 
 independent of $k$, such that
 \begin{equation} 
  u(r) \leq C \cl_{u}(2^{k} r), \qquad \forall r\geq 0.
     \end{equation}\notag
  \par\noindent
  (3) Let $u\in C_{+, \log}$ be increasing and 
	(log, $x^{k}$)-convex. Then for any $a>1$, we have
	\begin{equation} \label{eq:4-11}
	 \cl_{u}(r) \leq \sqrt{\ell_{u}(0){ea \over \log a}}\>
   	u\big(a2^{k+1}r\big)^{1/2}.  
\end{equation}
\end{fact}

\begin{fact} \label{fact:2-2}
If $u\in C_{+, 1/2}$ is (log, $x^{2}$)-convex, then the
Legendre transform $\,\ell_{u^{*}}$ of $u^{*}$ is given by
\begin{equation}
 \ell_{u^{*}} (t) = {e^{2t} \over \ell_{u}(t) t^{2t}},
 \qquad t\in [0, \infty).   \notag
\end{equation}
\end{fact}

\subsection{Complex Gaussian Space}\label{sec:2-2}
Let us start with taking a special choice of 
a Gel'fand triple:
\begin{equation*}
	\ce=\cs(\spr)\subset
	\ce_0=L^2(\spr,dt)
	\subset\ce^{*}=\cs^{*}(\spr)
\end{equation*}
just for convenience where $\cs$ is the Schwarz space 
of rapidly decreasing functions and $\cs^{*}$ is the 
space of tempered distributions.  
Consult \cite{hkps,kta,ktb,kuo96,ob}
for more general setting.
Let $A$ be  
a positive self-adjoint operator in $\ce_0$.  So there
exists an orthonormal basis 
$\{e_{j}\}_{j=0}^{\infty}\subset \ce$ for $\ce_0$ satisfying 
$Ae_{j}=\l_{j}e_{j}$.  
$|\cdot|_0$ denotes the norm of $\ce_0$.
For each $p\geq 0$ we define 
$|f|_p=|A^pf|_0$ and let 
$\ce_p=\{f\in \ce_0;\ |f|_p<\infty, \ p\geq 0\}$.
Note that $\ce_p$ is the completion of $\ce$ with 
respect to the norm $|\cdot|_p$.  Moreover, 
\begin{equation*}
	\r=\|A^{-1}\|_{OP},\quad
	\|i_{q,p}\|^2_{HS}=\sum_{j=0}^{\infty}\l^{-(q-p)}_j<\infty
\end{equation*}
for any $q>p\geq 0$.  Then the projective limit space $\ce$
of $\ce_p$ as $p\to \infty$ 
is a nuclear space and the dual space of $\ce$
is nothing but the inductive limit space $\ce'$.  
Hence we have the following 
continuous inclusions:
\begin{equation*}
	\ce\subset\ce_p\subset\ce_0\subset\ce'_p\subset\ce',
	\quad p\geq 0,
\end{equation*}
where the norm on $\ce'_p$ is given by 
$|f|_{-p}=|A^{-p}f|_0$.
Throughout this paper, we denote the complexification of a 
real space $X$ by $X_c$.
Let $\m'$ be the Gaussian measure on $\ce'$ with 
variance $1/2$, namely, a probability measure on $\ce'$ given 
by the characteristic function:
\begin{equation*}
	e^{-\frac{1}{4}|\x|^2_0}
	=\int_{\ce'}e^{i\la x,\x\ra}\m'(dx), 
	\qquad \x\in \ce.
\end{equation*}
Due to the topological isomorphism 
$\ce'_c\cong\ce'\times\ce'$, we can define a probability 
measure $\m_c=\m'\times\m'$ on $\ce'_c$.
The probability space $(\ce'_c,\m_c)$ is called the 
{\it complex Gaussian space}, see \cite{hida80}.
We denote by $L^2(\ce'_c,\m_c)$ the space of 
$\m_c$-square integrable functions on $\ce'_c$.
We should note that $L^2(\ce'_c,\m_c)\cong 
L^2(\ce',\m')\otimes L^2(\ce',\m')$.

Let $\m$ be the Gaussian measure on $\ce'$ with 
variance $1$ and $(\ce',\m)$ be the {\it real Gaussian space}.
The next Fact \ref{fact:2-3} has been obtained in \cite{akk7} for 
the Gel'fand triple $[\ce]_{u}\subset L^2(\ce',\m) \subset 
[\ce]_{u}^{*}$ associated with a growth function $u$.
This triple is refered to as 
{\it the CKS-space with a weight sequence 
$\a_{u}(n) = \big(\ell_{u}(n)n!)^{-1}$}.
For more precise discussion,  we will need the 
following conditions on $u \in C_{+, 1/2}$: 
\begin{itemize}
\item[(U0)] $\inf_{r\geq 0} u(r) = 1$.
\item[(U1)] $u$ is increasing and $u(0)=1$.
\item[(U2)] $\lim_{r\to\infty} r^{-1} \log u(r) < \infty$.
\item[(U3)] $u$ is (log, $x^{2}$)-convex.
\end{itemize}
Then we have
\begin{fact}\label{fact:2-3}
Suppose $u\in C_{+, 1/2}$ satisfies conditions (U0) (U2) 
(U3). Then the CKS-space with a weight sequence 
$\a_{u}(n)$ can be constructed.
Moreover, characterization theorems hold.
\end{fact}
\noindent
{\it Remark.}
(1) We refer the reader to the papers of Asai et al. 
\cite{akk6,akk7} for details.
We cite papers 
\cite{akk1,akk2,cks,ks92,ks93,kswy,kubo,kps,ps}
for characterizaion theorems on papticular cases.\\
(2) We should mention here that 
our formulation has some 
links with a recent work  
by Gannoun et al. \cite{ghor}.
So let us explain some of them.
Essential relationships are 
\begin{equation*}
u(r) = e^{2\theta(\sqrt{r})}, \quad 
u^*(r) = e^{2\theta^*(\sqrt{r})}
\end{equation*}
where $\t^{*}(s)=\sup_{t>0}\{st-\t(t)\}$
is adopted in \cite{ghor}.
In the following table
we give the correspondence
between our $U$-conditions and $\t$-conditions.
\begin{table}[h]
\begin{center}
	\begin{tabular}{|l|c|c|}
	\noalign{\hrule height0.8pt}
	\hfil  \ & $u$ &
	$\theta$ \\
	\hline
	$(U0)$ & $\displaystyle \inf_{r\geq 0} u(r) = 1$ 
	& $\displaystyle \inf_{r\geq 0} \theta(r) = 0$ \\
	$(U1)$ & $u$ is increasing and $u(0)=1$
	& $\theta$ is increasing and 
	$\theta(0) = 0$\\
	$(U2)$ & $\displaystyle \lim_{r\to\infty} \frac{\log u(r)}{r} < \infty$
	& $\displaystyle \lim_{r\to\infty} \frac{\theta(r)}{r^2} < \infty$\\
	$(U3)$ & $u$ is (log, $x^{2}$)-convex 
	& $\theta$ is convex\\
	\noalign{\hrule height0.8pt}
	\end{tabular}
\end{center}
\end{table}

\subsection{CKS-space over Complex Gaussian Space}
\label{sec:2-3}
Next let us consider the Gel'fand 
triple over the complex white noise space
$(\ce'_c,\m_c)$ for our purpose. 

$L^2(\ce'_c,\m_c)$-norm $\|\vf\|$ of $\vf$
is given by 
\begin{equation*}
	\|\vf\|^2=\sum_{l,m=0}^{\infty}l!m!|f_{l,m}|^2_0,
	\qquad f_{l,m}\in\ce^{\widehat{\otimes}(l+m)}_{0,c}.
\end{equation*}
In order to define norms in the spaces of test 
and generalized functions, we need a notation.
For $\k_{l,m}\in (\ce^{\otimes (l+m)}_c)^{*}_{symm}$,
we put 
\begin{equation*}
	|\k_{l,m}|_{p_1,p_2}=
	|(A^{p_1})^{\otimes l}
	\otimes(A^{p_2})^{\otimes m}\k_{l,m}|_0.
\end{equation*}
For $p_1,p_2\geq 0$ and given functions $u_1,u_2\in C_{+,1/2}$
satisfying conditions $(U0)(U2)(U3)$, define the
norm by 
\begin{equation}\label{eq:norm}
	\|\vf\|^2_{p_1,p_2}=\sum_{l,m=0}^{\infty}
	\frac{1}{\ell_{u_1}(l)\ell_{u_2}(m)}
	|f_{l,m}|^2_{p_1,p_2},
	\quad f_{l,m}\in \ce^{\widehat{\otimes}(l+m)}_c.
\end{equation}
Let 
	$[\ce_{p_1}]_{u_1}\otimes [\ce_{p_2}]_{u_2}
	=\bigl\{\vf\in L^2(\ce'_c,\m_c):\ \|\vf\|_{p_1,p_2}<\infty \bigr\}$.
Define the space $[\ce]_{u_1}\otimes [\ce]_{u_2}$
on $\ce'_c$ to be the projective limit of 
$[\ce_{p_1}]_{u_1}\otimes [\ce_{p_2}]_{u_2}$
as $p_1,p_2\to \infty$.
For abbreviation, we put 
$\cw_{u_1,u_2}=[\ce]_{u_1}\otimes [\ce]_{u_2}$ and 
$\cw^{*}_{u_1,u_2}=[\ce]^{*}_{u_1}\otimes [\ce]^{*}_{u_2}$
for the dual space.
For given functions $u_1,u_2\in C_{+,1/2}$ satisfying 
$(U0)(U2)(U3)$, we have the following continuous inclusions:
\begin{multline*}
	\cw_{u_1,u_2}\subset
	[\ce_{p_1}]_{u_1}\otimes [\ce_{p_2}]_{u_2}\subset
	L^2(\ce'_c,\m_c)\subset
	[\ce_{p_1}]^{*}_{u_1}\otimes [\ce_{p_2}]^{*}_{u_2}
	\subset 
	\cw^{*}_{u_1,u_2}.
\end{multline*}
where $[\ce_{p_1}]^{*}_{u_1}\otimes [\ce_{p_2}]^{*}_{u_2}$
is the dual space of $[\ce_{p_1}]_{u_1}\otimes [\ce_{p_2}]_{u_2}$.
In general, $u_1$ and $u_2$ are not necessarily 
the same functions.
A Gel'fand triple 
$\cw_{u_1,u_2}\subset 
L^2(\ce'_c,\m_c)\subset\cw^{*}_{u_1,u_2}$
is refered to as a {\it CKS-space with a weight sequence
$\a_{u_1}(l)\a_{u_2}(m)$}.
The bilinear form on $\cw^{*}_{u_1,u_2}\times \cw_{u_1,u_2}$
is denoted by $\lla\cdot,\cdot\rra_{c}$.
Then 
\begin{equation*}
	\lla \vF,\vf\rra_c=\sum_{n=0}^{\infty}
	l!m!\la F_{l,m},f_{l,m}\ra,
\end{equation*}
and it holds that 
\begin{equation*}
	|\lla \vF,\vf\rra_c|\leq \|\vF\|_{-p_1,-p_2}\|\vf\|_{p_1,p_2}
\end{equation*}
where
\begin{equation}\label{eq:-norm}
	\|\vF\|^2_{-p_1,-p_2}=\sum_{l,m=0}^{\infty}
	\frac{1}{\ell_{u^{*}_1}(l)\ell_{u^{*}_2}(m)}
	|F_{l,m}|^2_{-p_1,-p_2},
	\quad F_{l,m}\in (\ce^{\otimes (l+m)}_c)^{*}_{symm}.
\end{equation}

\subsection{Examples}\label{sec:2-4}
The combinations of two functions out of following examples
are applicable to our setting.
\begin{example}
Consider 
\begin{equation*}
	u(r)=u^{*}(r)=e^r.
\end{equation*}
Then it is obvious to check 
that consditions (U0) (U2) (U3) are satisfied.
This example produces to the {\it Hida-Kubo-Takenaka space}
over the real Gaussian space.
See \cite{hkps,kta,ktb,ob}.
\end{example}

\begin{example}
For $0\leq \b<1$, let $u$ be the function defined by
\begin{equation}
  u(r)= \exp\left[(1+\b)r^{1\over 1+\b}\right].  \notag
\end{equation}
It is easy to check that $u$ belongs to $C_{+, 1/2}$ and 
satisfies conditions (U0) (U2) (U3). By Example 4.3 in 
\cite{akk3}, the dual Legendre transform $u^{*}$ of $u$ 
is given by
\begin{equation}
 u^{*}(r) = \exp\left[(1-\b)r^{1\over 1-\b}\right]. \notag
\end{equation}
This example is for the construction of the {\it Kodratiev-Streit space}
over the real Gaussian space.
See \cite{ks92,ks93,kuo96}
\end{example}

\begin{example}
Consider the function $v(r) = \exp\big[e^{r}-1\big]$. 
Obviously, $v\in C_{+, 1/2}$. Let $u=v^{*}$ be the dual
Legendre transform of $v$. Then $u(0)=\sup_{s\geq 0}
v(s)^{-1}=1$ and it can be shown that $u$ belongs to
$C_{+, 1/2}$ and is an increasing (log, $x^{2}$)-convex 
function on $[0, \infty)$ (See \cite{akk6}). 
Hence $u\in C_{+, 1/2}$ 
satisfies conditions (U1) and (U3). It is shown in Example
4.4 in \cite{akk3} that $u$ is equivalent to the function
\begin{equation*} \label{eq:3-20}
 w(r) = \exp\left[2\sqrt{r\log\sqrt{r}}\,\right]. 
\end{equation*}
(``$u$ is equivalent to $v$" means that 
there exist constants $a_1,a_2,b_1,b_2>0$ 
and $r_0\in [0,\infty)$ such that 
$a_1b_1^nu(r)\leq v(r)\leq a_2b_2^nu(r)$ for all $r\geq r_0$.)
Obviously, $w$ satisfies condition (U2) and so $u$ also
satisfies condition (U2). On the other hand, we have 
the involution property $u^{*}=(v^{*})^{*}=v$. 
This example
can be applied to the Gel'fand triple $[\ce]_{u}\subset 
(L^{2}) \subset [\ce]_{u}^{*}$ for the following pair of
functions:
\begin{equation*}
 u^{*}(r) = \exp\big[e^{r}-1\big], \quad
 u(r) = (u^{*})^{*}.  \notag
\end{equation*}
In general, we can consider the follwing general pair of 
functions:
\begin{equation*}
	\exp_{k}(r)=\exp(\exp(\cdots(\exp(r)))),\quad
	w_{k}(r) = \exp\left[2\sqrt{r\log_{k-1}\sqrt{r}}\,\right]
\end{equation*}
We refer the reader to papers
\cite{asai1,akk1,akk2,akk3,akk4,akk6,akk7,cks,kks}.
\end{example}

\section{Characterizations of 
Product measures and Pseudo-positive Operators}
\label{sec:3}
We shall define another norm as follows.
Let $\cd_{p_1,p_2}$ for 
$p_1,p_2\geq 1$ be the space of all functions
$\varphi$ on $\ce^*_{c}\times\ce^*_{c}$ 
satisfying the following conditions:\\
(L1) $\varphi$ is an analytic function on 
$\ce^{*}_{p_1,c}\times\ce^{*}_{p_2,c}$. \\
(L2) There exists a nonnegative constant $C$ such that
\begin{equation*}
	|\varphi(x,y)|^2\leq Cu_1(|x|^2_{-p_1})u_2(|y|^2_{-p_2})
	\ \ \hbox{for any}\ 
	(x,y)\in \ce^{*}_{p_1,c}\times\ce^{*}_{p_2,c}.
\end{equation*}
For $\varphi\in \cd_{p_1,p_2}$, 
its norm is defined by
\begin{equation}\label{eq:Anorm}
	|\!|\!|\vf|\!|\!|_{p_1,p_2}
	:=\sup_{(x,y)\in \ce^{*}_{p_1,c}\times\ce^{*}_{p_2,c}}
	|\varphi(x,y)|u_1(|x|^2_{-p_1})^{-\frac{1}{2}}
	u_2(|y|^2_{-p_2})^{-\frac{1}{2}}.
\end{equation}
for a function $u\in C_{+,\log}$.
Define the space $\cd_{u_1,u_2}$ of
test functions on $\ce^{*}\times \ce^{*}$ 
to be the projective limit of $\cd_{p_1,p_2}$
as $p_1,p_2\to \infty$.
Let $\cd^{*}_{u_1,u_2}$ be 
the dual space of $\cd_{u_1,u_2}$.

\smallskip
\noindent
{\it Remark.}
This construction is motivated by 
Lee \cite{lee} and Asai et al. \cite{akk4,akk6,akk7}.
See also \cite{kon80,kuo96} and references cited therein.
Asai et al. \cite{akk4,akk6,akk7} and Gannoun et al. \cite{ghor}
have considered the case of $u_2\equiv 1$, independently.
In addition, Ouerdiane studied similar situations 
and the case $u_1(r^2)=u_2(r^2)=\exp(k^{-1}r^k)$ where $1\leq k\leq 2$.

\begin{lemma}\label{lem:3-1}
If $u_1, u_2\in C_{+,1/2}$ and an entire function $F(\x,\y)$
on $\ce_c\times\ce_c$ satisfies the growth condition
\begin{equation}\label{eq:UU}
	|F(\x,\y)|^2\leq C
	u^{*}_1(K_1|\x|^2_{p_1})u^{*}_2(K_2|\y|^2_{p_2})
\end{equation}
for a fixed positive $p_i\in \spr$, then for $q_i>p_i$ 
with $K_ie^2\|i_{q_i,p_i}\|^2_{HS}<1$, there exists a kernel 
$\k_{l,m}\in (\ce^{\otimes (l+m)}_c)^{*}_{symm}$
such that 
\begin{equation}
	F(\x,\y)=\sum_{l,m=0}^{\infty}
	\la \k_{l,m}, \y^{\otimes l}\otimes\x^{\otimes m}\ra
\end{equation}
and 
\begin{equation}
	|\k_{l,m}|^2_{-q_1,-q_2}
	\leq C^2(K_1e^2\|i_{q_1,p_1}\|^2_{HS})^l
	(K_2e^2\|i_{q_2,p_2}\|^2_{HS})^m
	\ell_{u^{*}_1}(l)\ell_{u^{*}_2}(m).
\end{equation}
\end{lemma}

\begin{proof}
Consider an entire function on $\spc^{m+l}$  
\begin{equation*}
	\psi=\psi(z_1,\cdots,z_m,w_1,\cdots,w_l)
	:=F(z_1\xi_1+\cdots +z_m\xi_m,w_1\eta_1+\cdots +w_l\eta_l).
\end{equation*}
Define an $(l+m)$-linear functional $V_{l,m}$ 
on $\ce_{c}\times \ce_{c}$
\begin{multline*}
	V_{l,m}(\xi_1,\cdots,\xi_m,\eta_1,\cdots.\eta_l)
	:=  \frac{1}{l!m!}
	\frac{\partial^{l+m}\psi}
	{\partial z_1\cdots\partial z_m\partial w_1\cdots\partial w_l}
	\Bigg |_{\scriptstyle z_1=\cdots=z_m=0 
	\atop \scriptstyle w_1=\cdots=w_l=0}  \\ 
	= \frac{1}{l!m!}\frac{1}{(2\pi)^{l+m}}
	\prod_{j=1}^{m}\int_{|z_j|=r_j}
	\frac{dz_j}{z^2_j}
	\prod_{k=1}^{l}\int_{|w_k|=r_k}
	\frac{dw_j}{w^2_k}
	\psi(z_1,\cdots,z_m,w_1,\cdots,w_l).
\end{multline*}
Taking $r=r_1|\xi_1|_{p_1}=\cdots =r_n|\xi_m|_{p_1}$ and 
$s=s_1|\eta_1|_{p_2}=\cdots =s_l|\eta_l|_{p_2}$
we get
\begin{equation*}
\begin{split}
	|V_{l,m}(\xi_1,\cdots,\xi_m,\eta_1,\cdots,\eta_l)|
	& \leq C\frac{1}{l!m!}
	\frac{u^{*}_1(K_1m^2r^2)^{1\over 2}}{r^m}
	\frac{u^{*}_2(K_2l^2s^2)^{1\over 2}}{s^l} \\
	& \hspace{2cm} \times |\xi_1|_{p_1} \cdots |\xi_m|_{p_1}
	|\eta_1|_{p_2} \cdots |\eta_l|_{p_2}
\end{split}
\end{equation*}
by \eqref{eq:UU}. 
Minimizing the right term, we have 
\begin{equation*}
\begin{split}
	|V_{l,m}(\xi_1,\cdots,\xi_m,\eta_1\cdots,\eta_l)|
	& \leq CK_1^{m\over 2}K_2^{l\over 2}{m^m \over m!}{l^l \over l!}
	\ell_{u^{*}_1}(m)^{1\over 2}
	\ell_{u^{*}_2}(l)^{1\over 2}  \\
	& \hspace{2cm} \times |\xi_1|_{p_1} \cdots |\xi_m|_{p_1}
	|\eta_1|_{p_2} \cdots |\eta_l|_{p_2}.
\end{split}
\end{equation*}
This shows that $V_{l,m}$ can be expressed in the form
\begin{equation*}
	V_{l,m}(\xi_1,\cdots,\xi_m,\eta_1,\cdots,\eta_l) 
	= \langle \kappa_{l,m}, 
	\eta_1\cdots\otimes\eta_l\otimes\zeta_1\otimes\cdots\otimes\zeta_m
	\rangle 
\end{equation*}
with 
	$\kappa_{l,m} \in (\ce^{\otimes(l+m)}_{c})^*_{symm}$, 
\begin{equation*}
	\vert \kappa_{l,m}\vert^2_{-(p+q)}
	\leq C^2
	(K_1e^2\| i_{q_1,p_1}\|^2_{HS})^{l}
	(K_2e^2\| i_{q_2,p_2}\|^2_{HS})^{m}
	\ell_{u^{*}_1}(l)\ell_{u^{*}_2}(m)
\end{equation*}
with finite 
Hilbert-Schmidt norm $K_ie^2\|i_{q_i,p_i}\|^2_{H.S.} < 1$ 
for any $q_i>p_i$.  
Therefore we derive 
\begin{equation*}
	F(\xi,\eta) = \sum_{l,m=0}^\infty \langle \kappa_{l,m}, 
	\eta^{\widehat{\otimes}l}\otimes\zeta^{\widehat{\otimes}m}\rangle.
\end{equation*}
\end{proof}

Similarly, we have 
\begin{lemma}\label{lem:3-2}
If $u_1, u_2\in C_{+,\log}$ and an entire function $F(\x,\y)$
on $\ce'_c\times\ce'_c$ satisfies the growth condition
\begin{equation}\label{eq:uu}
	|F(\x,\y)|^2\leq Cu_1(K_1|\x|^2_{-p_1})u_2(K_2|\y|^2_{-p_2})
\end{equation} 
for any $K_i,p_i\geq 0 \ (i=1,2)$ and some $C>0$, then 
there exists a kernel 
$\k_{l,m}\in \ce^{\widehat{\otimes}(l+m)}_c$ such that 
\begin{equation}
	F(\x,\y)=\sum_{l,m=0}^{\infty}
	\la \k_{l,m}, \y^{\otimes l}\otimes \x^{\otimes m}\ra
\end{equation}
and 
\begin{equation}
	|\k_{l,m}|^2_{q_1,q_2}
	\leq C^2(K_1e^2\|i_{p_1,q_1}\|^2_{HS})^l
	(K_2e^2\|i_{p_2,q_2}\|^2_{HS})^m
	\ell_{u_1}(l)\ell_{u_2}(m).
\end{equation}
for any $q_i<p_i \ (i=1,2)$ satisfying 
$K_ie^2\|i_{p_i,q_i}\|^2<1$.
\end{lemma}

\begin{proof}
We can prove this Lemma with modifications of the 
proof of the previous Lemma.  Therefore, we omit the proof.
\end{proof}

\noindent
{\it Remark.}
Lemma \ref{lem:3-2} will be expected to characterize the 
continuous linear operator from $[\ce]^{*}_{u_2}$
into $[\ce]_{u_1}$ and expand it in terms of integral kernel 
operators \cite{ob,ob99}.  To our best knowledge, 
such consideratoins have not been done in any literature.

For $\vF\in \cw_{u_1,u_2}$, its {\it multiple S-transform}
$S_m\vF$ is defined to be the function
\begin{equation}
	(S_m\vF)(\x,\y)
	=\lla \vF, 
	e^{\sqrt{2}\la x,\x\ra-{1\over 2}\la\x,\x\ra}\otimes
	e^{\sqrt{2}\la y,\y\ra-{1\over 2}\la\y,\y\ra}
	\rra_c, \quad \x,\y\in\ce_c.
\end{equation}
For $\vf\in L^2(\ce'_c,\m')$, we have an integral representation 
of the multiple S-transform as
\begin{equation}
	S_m\varphi(\x,\y)
	=\int_{\ce^{*}}\int_{\ce^{*}}
	\varphi(x+\sqrt{2}\x,y+\sqrt{2}\y)
	\m'(dx)\m'(dy).
\end{equation}
Note that the multiple S-transform is 
nothing but a {\it symbol of operators}
frequently used in white noise operator theory
\cite{ob,ob99}.  Lemma \ref{lem:3-1}, Equation \eqref{eq:UU}
and $S_m$
(Lemma \ref{lem:3-2}, Equation \eqref{eq:uu} and $S_m$)
give us the characterization theorem for
 $\cw^{*}_{u_1,u_2}$ ($\cw_{u_1,u_2}$),
respectively, which generalize recent results
in \cite{akk7}.

\begin{proposition}\label{prop:3-3}
Suppose $u_1, u_2\in C_{+,1/2}$ satisfy 
$(U0)$$(U2)$$(U3)$.  Then the families of norms 
$\{|\!|\!|\cdot|\!|\!|_{p_1,p_2};\ p_1, p_2\geq 0\}$ 
and 
$\{\|\cdot\|_{p_1,p_2};\ p_1, p_2\geq 0\}$ 
are equivalent.
\end{proposition}

\begin{proof}
First, we will show that for any $p_i\geq 1$, $i=1,2$,
there exist $C\geq 0$ and $q_i>p_i$ such that
$|\!|\!|\varphi|\!|\!|_{p_1,p_2}
\leq C\|\varphi\|_{q_1,q_2}$.
Let $p_i\geq 1$, $i=1,2$, be given. 
Since it has been proved \cite{lee} (see also \cite{akk7,kuo96})
that every test function in $[\ce]_u$ 
has an analytic extention,   
there exist $C\geq 0$ and $q_i\geq p_i$ such that
for any $\varphi\in [\ce_{q_1}]_{u_1}\otimes [\ce_{q_2}]_{u_2}$
\begin{equation}\label{eq:phiestimate}
	|\varphi(x,y)|
	\leq C\|\varphi\|_{q_1,q_2}u_1(|x|^2_{-p_1})^{{1 \over 2}}
	u_2(|y|^2_{-p_2})^{{1 \over 2}}
\end{equation}
for any $(x,y)\in \ce^{*}_{p_1,c}\times\ce^{*}_{p_2,c}$.
Hence it is derived by \eqref{eq:Anorm} 
and \eqref{eq:phiestimate} that
\begin{align}\label{eq:AE}
	|\!|\!|\varphi|\!|\!|_{p_1,p_2}
	& =\sup_{(x,y)\in\ce^{*}_{p_1,c}\times\ce^{*}_{p_2,c} }
	|\varphi(x,y)|u(|x|^2_{-p_1})^{-\frac{1}{2}}
	u(|y|^2_{-p_2})^{-\frac{1}{2}} \notag \\
	& \leq C\|\varphi\|_{q_1,q_2}.
\end{align}

To prove the converse, 
The multiple S-transform of $\varphi$ is given by
$$
	F(\x,\y):=S_m\varphi(\x,\y)
	=\int_{\ce^{*}}\int_{\ce^{*}}
	\varphi(x+\sqrt{2}\x,y+\sqrt{2}\y)
	\m'(dx)\m'(dy).
$$
Then observe that for $q_i\geq 1$
\begin{align*}
	| & F(\x, \y)| 
	\leq \int_{\ce^{*}}\int_{\ce^{*}}
	|\varphi(x+\sqrt{2}\x,y+\sqrt{2}\y)|
	\m'(dx)\m'(dx) \\
	& = \int_{\ce^{*}}\int_{\ce^{*}}
	|\varphi(x+\sqrt{2}\x,y+\sqrt{2}\y)|
	u\bigl(|x+\sqrt{2}\x|^2_{-q_1}\bigr)^{-\frac{1}{2}}
	u\bigl(|y+\sqrt{2}\y|^2_{-q_2}\bigr)^{-\frac{1}{2}} \\
	& \hspace*{3cm}
	\times u\bigl(|x+\sqrt{2}\x|^2_{-q_1}\bigr)^{\frac{1}{2}}
	u\bigl(|y+\sqrt{2}\y|^2_{-q_2}\bigr)^{\frac{1}{2}}
	\m'(dx)\m'(dy) \\
	& \leq |\!|\!|\varphi|\!|\!|_{q_1,q_2}
	\int_{\ce^{*}}\int_{\ce^{*}}
	u\bigl(|x+\sqrt{2}\x|^2_{-q_1}\bigr)^{\frac{1}{2}}
	u\bigl(|y+\sqrt{2}\y|^2_{-q_2}\bigr)^{\frac{1}{2}}
	\m'(dx)\m'(dy). 
\end{align*}
By the condition $(U0)$, $u^{\frac{1}{2}}(r)\leq u(r)$ 
for all $r\geq 0$.
Therefore, 
\begin{equation*}
	|F(\x, \y)|
	\leq \int_{\ce^{*}}\int_{\ce^{*}}
	u\bigl(|x+\sqrt{2}\x|^2_{-q_1}\bigr)
	u\bigl(|y+\sqrt{2}\y|^2_{-q_2}\bigr)
	\m'(dx)\m'(dy)
\end{equation*}
By the condition $(U3)$, we have
\begin{align*}
	u(|x+\sqrt{2}\x|^2_{-q})
	& \leq u\bigl((|x|_{-q}
	+|\sqrt{2}\x|_{-q})^2\bigr) \\
  	& \leq u(4|x|^2_{-q})^{{1\over 2}}
  	u(8|\x|^2_{-q})^{{1\over 2}}.
\end{align*}
Thus, it is easy to get
\begin{align}\label{eq:uugrowth1}
	|F(\x,\y)|
	& \leq L|\!|\!|\varphi|\!|\!|_{q_1,q_2}
	u_1(8|\x|^2_{-q_1})^{{1\over 2}}
	u_2(8|\y|^2_{-q_2})^{{1\over 2}}
\end{align}
where 
\begin{equation*}
	L=\int_{\ce^{*}}\int_{\ce^{*}}
	u(4|x|^2_{-q_1})^{{1\over 2}}
	u(4|y|^2_{-q_2})^{{1\over 2}}
	\m'(dx)\m'(dy)<\infty.
\end{equation*}
(Note that finiteness concerning $L$ can be shown easily
by simple estimation 
and the Fernique theorem \cite{fer, kuo96}.)
Then applying Lemma \ref{lem:3-2} with \eqref{eq:uugrowth1},
we have 
\begin{align}
	\|\varphi\|^2_{p_1,p_2}
	& \leq L^2(1-8e^2\|i_{q_1,p_1}\|^2_{HS})^{-1}
	(1-8e^2\|i_{q_2,p_2}\|^2_{HS})^{-1}
	|\!|\!|\varphi|\!|\!|^2_{q_1,q_2}. \label{eq:EA}
\end{align}
We complete the proof.
\end{proof}

\begin{definition}\label{def:positive}
A generalized function 
$\vF\in \cw^{*}_{u_1,u_2}$
is called {\it positive} if 
$\lla \vF, \vf\rra\geq 0$
for all nonnegative test functions 
$\vf\in \cd_{u_1,u_2}$. 
\end{definition}
\noindent
{\it Remark.}
Positivity of generalized functions in 
white noise context has been studied by Yokoi \cite{yokoi90}.\\

\smallskip
It is possible to give an alternative definition
to Definition \ref{def:positive}
as follows by the kernel theorem \cite{treves}.
\begin{definition}
An operator $\varXi\in\cl([\ce]_{u_1},[\ce]^{*}_{u_2})$ 
is called {\it positive} in the sense of 
distributions if $\lla \varXi\vf_1, \vf_2\rra_c\geq 0$
for all nonnegative test functions 
$\vf_i\in[\ce]_{u_i}, \ (i=1,2)$.
We call such an operator {\it pseudo-positive} operator
\end{definition}
\noindent
{\it Notation.}
For locally convex spaces $X, Y$, let $\cl(X,Y)$
denote the space of all continuous operators from 
$X$ into $Y$ equipped with the topology 
of uniform convergence on every bounded subset.

\begin{theorem}
Suppose $u_1,u_2\in C_{+,1/2}$ satisfy
$(U0)(U2)(U3)$.  A measure $\n_1\times\n_2$ on 
$\ce'\times\ce'$ is a positive product Radon measure 
inducing a positive generalized function 
$\vF_{\n_1\times\n_2}\in \cw^{*}_{u_1,u_2}$
if and only if $\n_1\times\n_2$ is supported in 
$\ce^{*}_{p_1}\times\ce^{*}_{p_2}$ for some $p_1,p_2\geq 1$
and 
\begin{equation}\label{eq:intcond}
	\int_{\ce^{*}_{p_1}}\int_{\ce^{*}_{p_2}}
	u_1(|x|^2_{-p_1})^{\frac{1}{2}}u_2(|y|^2_{-p_2})^{\frac{1}{2}}
	\n_1(dx)\n_2(dy)<\infty.
\end{equation}
\end{theorem}

\begin{proof}
First we shall prove sufficiency.
Suppose that $\n_1\times\n_2$ is supported in 
$\ce'_{p_1}\times\ce'_{p_2}$ for some $p_1, p_2\geq 0$
and Equation \eqref{eq:intcond} holds.  Then for any 
$\vf(x,y)\in \cw_{u_1,u_2}$, 
\begin{align}
	& \int_{\ce'_{p_1}}\int_{\ce'_{p_2}}
	|\vf(x,y)|\n_1(dx)\n_2(dy)\notag \\
	& = \int_{\ce'_{p_1}}\int_{\ce'_{p_2}}
	|\vf(x,y)|u_1(|x|^2_{-p_1})^{-\frac{1}{2}}
	u_2(|y|^2_{-p_2})^{-\frac{1}{2}} \notag\\
	& \hspace{3cm}\times u_1(|x|^2_{-p_1})^{\frac{1}{2}}
	u_2(|y|^2_{-p_2})^{\frac{1}{2}} 
	\n_1(dx)\n_2(dy)\notag \\
	& \leq |\!|\!|\vf|\!|\!|_{p_1,p_2}
	\int_{\ce'_{p_1}}\int_{\ce'_{p_2}}
	u_1(|x|^2_{-p_1})^{\frac{1}{2}}
	u_2(|y|^2_{-p_2})^{\frac{1}{2}} 
	\n_1(dx)\n_2(dy)
\end{align}
With the help of Proposition \ref{prop:3-3}, 
$\cw_{u_1,u_2}
\subset L^1(\ce'_c,\n_1\times\n_2)$ and 
\begin{equation}
	\vf\longmapsto \int_{\ce'_{p_1}}\int_{\ce'_{p_2}}
	\vf(x,y)\n_1(dx)\n_2(dy)
\end{equation} 
is a continuous linear functional on 
$\cw_{u_1,u_2}$.
Therefore, $\n_1\times\n_2$ is a positive product Radon measure
which induces a positive generalized function 
$\vF_{\n_1\times\n_2}$ in $\cw^{*}_{u_1,u_2}$.

Conversely, suppose that $\n_1\times\n_2$ is a positive product 
Radon measure.  Then for all 
$\vf\in \cw_{u_1,u_2}$,
\begin{equation}
	\lla \vF_{\n_1\times\n_2}, \vf\rra
	=\int_{\ce'}\int_{\ce'}\vf(x,y)\n_1(dx)\n_2(dy).
\end{equation}
is a continuous linear functional with respect to 
$\{|\!|\!|\cdot|\!|\!|^2_{q_1,q_2}; \ p_1,p_2\geq 0\}$
by Proposition \ref{prop:3-3}.
Thus there exist constants $K,q_1,q_2\geq 0$ such that
for all $\vf\in \cd_{u_1,u_2}$
\begin{equation}\label{eq:functional}
	\bigl|\lla \vF_{\n_1\times \n_2}, \vf\rra\bigr|
	\leq K|\!|\!|\vf|\!|\!|_{q_1,q_2}.
\end{equation}

Let us define an analytic function $\om$ on 
$\ce'_{q_1c}\times\ce'_{q_2,c}$ by 
\begin{equation}
	\om(x,y)=\cl_{u_1}(2^{-4}\la x,x\ra_{-q_1})
	\cl_{u_2}(2^{-4}\la y,y\ra_{-q_2}), \
	(x,y)\in \ce'_{q_1c}\times\ce'_{q_2,c}
\end{equation}
where $\la\cdot,\cdot\ra_{-q_i}$ is the bilinear 
pairing on $\ce'_{q_i,c},\ (i=1,2)$.
On the other hand, Fact \ref{fact:2-1} (3) implies that 
\begin{align}
	|\om(x,y)|
	& \leq \cl_{u_1}(2^{-4}|x|^2_{-q_1})
	\cl_{u_2}(2^{-4}|x|^2_{-q_2})\notag \\
	& \leq \frac{2e}{\log 2}
	u_1(|x|^2_{-q_1})^{1\over 2}u_2(|x|^2_{-q_2})^{1\over 2}.
\end{align}
(It is easy to find an increasing function $v$ equivalent to a 
function $u$.)
This implies that $\om\in \cd_{p_1,p_2}$.
Thus, from Equation \eqref{eq:functional} with $\vf=\om$ we 
obtain that 
\begin{equation}\label{eq:omega}
	\biggl|\int_{\ce'}\om(x,y)\n_1(dx)\n_2(dy)\biggr|
	\leq K\frac{2e}{\log 2}.
\end{equation}
Due to Equation \eqref{eq:omega}, we have 
\begin{equation*}
	\int_{\ce'}\om(x,y)\n_1(dx)\n_2(dy)<\infty.
\end{equation*}
However Fact \ref{fact:2-1} (2) says that 
$u_1(r_1)u_2(r_2)\leq Cw(4r_1,4r_2)$.  Therefore,
\begin{equation*}
	\int_{\ce'}u_1(2^{-6}|x|^2_{-q_1})u_2(2^{-6}|y|^2_{-q_2})
	\n_1(dx)\n_2(dy)<\infty.
\end{equation*}
By choosing an appropriate $p_i>q_i \ (i=1,2)$ satisfying 
$\r^{2(p_i-q_2)}\leq 2^{-6}$, so that 
$|\cdot|^2_{-p_i}\leq 2^{-6}|\cdot|^2_{-q_i}\ (i=1,2)$.
Therefore we get the assertion.
\end{proof}

\begin{theorem}
Suppose $u_1,u_2\in C_{+,1/2}$ satisfy 
$(U0)(U2)(U3)$.  A measure $\n_1\times\n_2$ on 
$\ce'\times\ce'$ is a positive product Radon measure 
inducing a pseudo-positive operator 
$\vX\in\cl([\ce]_{u_1},[\ce]^{*}_{u_2})$
if and only if $\n_1\times\n_2$ is supported in 
$\ce^{*}_{p_1}\times\ce^{*}_{p_2}$ for some $p_1,p_2\geq 1$
and 
\begin{equation}
	\int_{\ce^{*}_{p_1}}\int_{\ce^{*}_{p_2}}
	u_1(|x|^2_{-p_1})^{\frac{1}{2}}u_2(|y|^2_{-p_2})^{\frac{1}{2}}
	\n_1(dx)\n_2(dy)<\infty.
\end{equation}
\end{theorem}


\end{document}